\documentclass[11pt]{article}
\usepackage{amsfonts,amsmath,amssymb}
\usepackage[dvips]{graphicx}
\usepackage[dvips]{epsfig}
\usepackage{hangcaption}
\usepackage[T1]{fontenc}

\usepackage{vmargin} 

\newcommand{\dd}{\delta}
\newcommand{\e}{\epsilon}
\newcommand{\alp}{\alpha}
\newcommand{\lab}{\lambda}
\newcommand{\gam}{\gamma}
\newcommand{\sig}{\sigma}
\newcommand{\tht}{\theta}

\newcommand{\sli}{\sum\limits}
\newcommand{\sliin}{\sum\limits_{i=1}^n}

\newcommand{\proliin}{\prod\limits_{i=1}^n}

\newcommand{\ili}{\int\limits}

\newcommand{\bculi}{\bigcup\limits}

\newcommand{\limn}{\lim_{n\rightarrow\infty}\;}
\newcommand{\limk}{\lim_{k\rightarrow\infty}\;}
\newcommand{\lsn}{\limsup_{n\rightarrow\infty}}
\newcommand{\lin}{\liminf_{n\rightarrow\infty}}
\newcommand{\lsk}{\limsup_{k\rightarrow\infty}}

\newcommand{\rar}{\rightarrow}

\newcommand{\kif}{k\rightarrow\infty}
\newcommand{\nif}{n\rightarrow\infty}

\newcommand{\aoo}{\Big\{}
\newcommand{\aff}{\Big\}}
\newcommand{\coo}{\Big [}
\newcommand{\cff}{\Big]}
\newcommand{\poo}{\Big (}
\newcommand{\pff}{\Big)}
\newcommand{\po}{\big (}
\newcommand{\pf}{\big)}
\newcommand{\ao}{\big \{}
\newcommand{\af}{\big \}}

\newcommand{\pooo}{\bigg (}
\newcommand{\pfff}{\bigg)}
\newcommand{\aooo}{\bigg \{}
\newcommand{\afff}{\bigg \}}
\newcommand{\cooo}{\bigg [}
\newcommand{\cfff}{\bigg ]}
\newcommand{\poooo}{\Bigg (}
\newcommand{\pffff}{\Bigg)}

\newcommand{\coooo}{\Bigg [}
\newcommand{\cffff}{\Bigg ]}
\newcommand{\mmi}{\mid\mid}

\newcommand{\Mid}{\Big |}

\newcommand{\nk}{n_k}
\newcommand{\nkm}{n_{k-1}}

\newcommand{\hk}{h_{n_k}}

\newcommand{\EEE}{\mathbb{E}}

\newcommand{\PPP}{ \mathbb{P}}

\newcommand{\RRR}{\mathbb{R}}

\newcommand{\QQQ}{\mathbb{Q}}
\newcommand{\AAA}{\mathcal{A}}
\newcommand{\BB}{\mathcal{B}}
\newcommand{\CC}{\mathcal{C}}
\newcommand{\DD}{\mathcal{D}}
\newcommand{\FF}{\mathcal{F}}
\newcommand{\GG}{\mathcal{G}}
\newcommand{\HH}{\mathcal{H}}
\newcommand{\KK}{\mathcal{K}}
\newcommand{\LL}{\mathcal{L}}
\newcommand{\MM}{\mathcal{M}}
\newcommand{\NN}{\mathcal{N}}

\newcommand{\RR}{\mathcal{R}}

\newcommand{\TT}{\mathcal{T}}

\newcommand{\ovg}{\overline{g}}
\newcommand{\ovh}{\overline{h}}
\newcommand{\ovr}{\overline{r}}

\newcommand{\wt}{\widetilde}

\newcommand{\wtd}{\widetilde{d}}

\newcommand{\wth}{\mathfrak{h}}

\newcommand{\wtz}{\widetilde{z}}

\newcommand{\wtW}{\widetilde{W}}
\newcommand{\wtY}{\widetilde{Y}}

\newcommand{\beq}{\begin{equation}}
\newcommand{\eeq}{\end{equation}}
\newcommand{\lb}{\newline}
\newcommand{\indep}{{\bot}\kern-0.9em{\bot}} 

\newcommand{\mcal}{\mathcal}

\newcommand{\wap}{(\Omega,\mathcal{A},\rm I\kern-2pt P)} 
\newcommand{\vk}{\vskip10pt}

\newcommand{\nono}{\nonumber}

\newcommand{\Var}{\mathrm{Var}}

\newcommand{\suite}{_{n\geq 1}}

\newcommand{\mh}{\mathfrak{h}}

\newcommand{\zklj}{\mathbf{z}}
\newcommand{\hkl}{h_{n_k,l}}
\newtheorem{theo}{Theorem}

\newtheorem{lem}{Lemma}

\newtheorem{coro}{Corollary}

\newtheorem{rem}{Remark}
\newtheorem{fact}{Fact}

\begin{document}
\captionwidth= 12cm
\title{Uniform in bandwidth exact rates \\ for a class of kernel estimators}
\author{{\large Davit} {\Large V}{\normalsize ARRON \footnote{Laboratoire de Mathématiques Pures et Appliquées, Université de Franche-Comté, 16 route de gray, 25000 Besançon, France, E-mail: dvarron@univ-fcomte.fr.  Financial support from IAP research network P6/03 of the Belgian Government (Belgian Science Policy).}}
\and
{\large Ingrid} {\Large V}{\normalsize AN} {\Large K}{\normalsize EILEGOM
\footnote{Institute of Statistics, Universit\'e catholique de Louvain, Voie du Roman Pays 20, 1348 Louvain-la-Neuve, Belgium, E-mail: ingrid.vankeilegom@uclouvain.be.  Financial support from IAP research network P6/03 of the Belgian Government (Belgian Science Policy), and from the European Research Council under the European Community's Seventh Framework Programme (FP7/2007-2013) / ERC Grant agreement No. 203650 is gratefully acknowledged.}}}

\date{\today}

\maketitle

\def\baselinestretch{1.1}

\begin{abstract}
Given an i.i.d sample $(Y_i,Z_i)$, taking values in
$\RRR^{d'}\times \RRR^d$, we consider a collection Nadarya-Watson kernel estimators of the conditional expectations $\EEE( <c_g(z),g(Y)>+d_g(z)\mid Z=z)$, where $z$ belongs to a
compact set $H\subset \RRR^d$, $g$ a Borel function on
$\RRR^{d'}$ and $c_g(\cdot),d_g(\cdot)$ are continuous functions
on $\RRR^d$. Given two bandwidth sequences $h_n<\wth_n$ fulfilling
mild conditions, we obtain an exact and explicit almost sure limit bounds for the deviations of these estimators around their expectations, uniformly in $g\in\GG,\;z\in H$ and $h_n\le h\le \wth_n$ under mild conditions on
the density $f_Z$, the class $\GG$, the kernel $K$ and the
functions $c_g(\cdot),d_g(\cdot)$. We apply this result to prove
that smoothed empirical likelihood can be used to build confidence
intervals for conditional probabilities $\PPP( Y\in C\mid
Z=z)$, that hold uniformly in $z\in H,\; C\in \CC,\; h\in
[h_n,\wth_n]$. Here $\CC$ is a Vapnik-Chervonenkis class of sets.\\

\end{abstract}

\def\baselinestretch{1.2}

\noindent {\large Key Words:} Local empirical processes, empirical likelihood, kernel smoothing,
uniform in bandwidth consistency.

\newpage
\normalsize

\section{Introduction and statement of the main results}
Consider an i.i.d sample $(Y_i,Z_i)_{i=1,\ldots,n}$ taking values
in $\RRR^{d'}\times \RRR^d$, with the same distribution as a
vector $(Y,Z)$, and write $<\cdot,\cdot>$ for the usual inner
product. In this paper, we investigate the limit behaviour of
quantities of the following form (assuming that this expression is
meaningful):
\begin{align}
 \nono W_n(g,h,z):=& f_Z(z)^{-1/2}\sliin \cooo \po
<c_g(z),g(Y_i)>+d_g(z)\pf K\poo \frac{Z_i-z}{h}\pff \\
 &\;- \EEE \pooo\po <c_g(z),g(Y_i)>+d_g(z)\pf K\poo
\frac{Z_i-z}{h}\pff\pfff\cfff.\label{Wnghz}\end{align} Here, $K$
denotes a kernel, $h>0$ is a smoothing parameter, $g$ is a Borel
function from $\RRR^{d'}$ to $\RRR^k$ and $f_Z$ is (a version) of
the density of $Z$. Given a class of functions $\GG$ satisfying
some Vapnik-Chervonenkis type conditions (see conditions (HG1)
below), and given a compact set $H$, Einmahl and Mason
(2000) showed that somewhat recent tools in empirical
processes theory could be used efficiently to provide exact rates
of convergence of
$$\sup\ao\mid W_n(g,h_n,z)\mid,\;g\in \GG,\;z\in H\af,$$
along a bandwidth sequence $h_n$ fulfilling some mild conditions
(see condition $(HV)$ in the sequel). The exact content of their
result is written in Theorem \ref{EM} below. The contribution of
the present paper is twofold. As a first contribution, we provide
an extension of the result of Einmahl and Mason, by enriching
Theorem \ref{EM} with a uniformity in the bandwidth $h$, when $h$
is allowed to vary into an interval $[h_n,\wth_n]$, with $h_n$ and
$\wth_n$ fulfilling conditions of Theorem \ref{EM}. This extension
is stated in Section \ref{emplus} (Theorem \ref{T1}), and is
proved in Section \ref{PT1}. As a second contribution (Theorem
\ref{T2}), we apply our Theorem \ref{T1} to establish confidence
intervals for quantities of the form
$$\PPP\poo Y\in C\mid Z=z\pff,\;C\in \CC,\;z\in H,$$
by empirical likelihood techniques. Indeed, we prove that these
confidence intervals can be built to hold uniformly in $z\in
H,\;C\in \CC$ and $h\in [h_n,\wth_n]$, under conditions that are
very similar to those of Theorem \ref{T1}. This result is stated
in Section \ref{elplus} and is proved in Section \ref{PT2}.
\subsection{A result of Einmahl and Mason}
As our first result is an extension of Theorem 1 in
Einmahl and Mason (2000) we have to first introduce the notations and
assumptions they made in their article. Consider a compact set
$H\subset \RRR^d$ with nonempty interior. We shall make the
following assumption on the law of $(Y,Z)$.
\begin{tabbing} $\hskip5pt$ \= $(H f)$
$\hskip5pt$ \= $(Y,Z)$ has a density $f_{Y,Z}$ that is continuous in $x$ on $\RRR^{d'}\times O'$, where $O'\subset \RRR^d$\\
\>\>  is open and where $H\subset O'$. \\
\>\>Moreover $f_Z$ is continuous and bounded away from zero and
infinity on $O'$.
\end{tabbing}
From now on, $O$ will denote an open set fulfilling $ H\subsetneq
O\subsetneq O'.$ Now consider a class $\GG$ of functions from
$\RRR^{d'}$ to $\RRR^k$. For $l=1,\ldots,k$, write
$\GG_l:=\Pi_l(\GG)$, where $\Pi_l(x_1,\ldots,x_l,\ldots,x_k):=x_l$
for $(x_1,\ldots,x_k)\in \RRR^k$.
\begin{tabbing}
$\hskip5pt$ \= $(HG)$$\hskip5pt$ \=  Each class $\GG_l$ is a
pointwise separable VC subgraph class and has a finite valued\\
\>\> measurable envelope function $G_l$ satisfying, for some $p\in (2,\infty]$:\\
\>\>$ \alp:=\max_{l=1,\ldots,k}\;\sup_{z\in O}\mmi
G_l(\cdot)\mmi_{\LL_{Y\mid Z=z},p}<\infty,$
\end{tabbing}
where $\mmi G_l(\cdot)\mmi_{\LL_{Y\mid Z=z},p}$ is the $L^p$-norm
of $G_l$ under the distribution of $Y\Mid Z=z$. For a definition
of a pointwise separable VC subgraph class we refer to
Van de Vaart and Wellner (1996, p. 110 and 141). Now, for any $g\in \GG$, consider
a pair of functions $(c_g(\cdot),d_g(\cdot))$, where $c_g$ maps
$\RRR^d$ to $\RRR^k$ and $d_g$ maps $\RRR^d$ to $\RRR$, and assume
that
\begin{tabbing}
$\hskip5pt$ \= $(H\CC)$$\hskip5pt$ \= The classes of functions
$\DD_1:=\{c_g,\;g\in \GG\}\;\text{and}\;\DD_2:=\{d_g,\;g\in \GG\}$
are uniformly\\
 \>\>bounded and uniformly equicontinuous on $O$.
\end{tabbing}
We now formulate our
assumptions on the Kernel $K$, with the following definition. \beq
\KK:=\aoo K\po \lab\cdot-z\pf,\;\lab>0,\;z\in
\RRR^d\aff.\label{KK}\eeq
 \begin{tabbing} $\hskip5pt$ \= $(HK1)$
$\hskip5pt$ \= $K$ has bounded variation and the class $\KK$ is VC subgraph. \\
\> $(HK2)$\> $K(s)=0$ when $s\notin [-1/2,1/2]^d$.\\
\> $(HK3)$ \> $\ili_{\RRR^d} K(s)ds=1$.
\end{tabbing}
Note that $(HK1)$ is fulfilled for a quite large class of kernels
(see, e.g., Mason (2004), Example F.1). In
Einmahl and Mason (2000), the authors have studied the almost sure
asymptotic behaviour of
$$\sup\ao\mid W_n(g,h_n,z)\mid,\;g\in \GG,\;z\in H\af$$ (recall (\ref{Wnghz})), along a bandwidth
sequence $(h_n)\suite$ that satisfies the following conditions
(here we write $\log_2 n:= \log \log (n\vee 3)$) :
\begin{tabbing}
$\hskip5pt$ \= $(HV)$$\hskip5pt$ \= $h_n\downarrow
0,\;n{h_n}^d\uparrow\infty,\;\log(1/h_n)/\log_2 n\rar\infty$,\;
$h_n^d \po n/\log(1/h_n)\pf^{1-2/p}\rar \infty$,
\end{tabbing}
where $p$ is as in condition $(HG)$. We also set
\begin{align}
\Delta^2(g,z):=&\EEE\poo \po <c_g(z),g(Y)>+d_g(z)\pf^2\Mid Z=z\pff,\;z\in \RRR^d,\;g\in \GG,\label{siggz}\\
\Delta^2(g):=&\sup_{z\in H} \Delta^2(g,z),\;g\in \GG\label{sigg}\\
\Delta^2(\GG):=&\sup_{g\in \GG}\Delta^2(g)\label{sigGG}.
\end{align}
Given a measurable space $(\chi,\TT)$, a measure $Q$ and a Borel
function $\psi: \;\chi\mapsto \RRR$, we write \beq {\mmi
\psi\mmi^p_{Q,p}}=\ili_{\chi} \mid\psi^p\mid dQ \label{Lp}.\eeq
Under the above mentioned assumptions, Einmahl and Mason have
proved the following theorem, $\lab$ denoting the Lebesgue
measure.
\begin{theo}[Einmahl, Mason, 2000]\label{EM}
Under assumptions $(HG)$, $(H \mcal{C})$, $(H f)$, $(HK1)-(HK3)$
and $(HV)$, we have almost surely \beq \limn \sup_{z\in H,\; g\in
\GG} \frac{\mid W_n(g,h_n,z)\mid}{\sqrt{2 n h^d_n
\log(h^{-d}_n)}}=\Delta(\GG)\mmi K\mmi_{\lab,2}.\label{ben}\eeq
\end{theo}
We point out that (\ref{ben}) is slightly stronger than Theorem 1
of Einmahl and Mason (2000) , as $f_Z(z)^{-1/2}$ appears
in our definition of $W_n(g,h,z)$ which is not the case in their
paper. However, (\ref{ben}) is a consequence of their Theorem 1,
as $f_Z^{-1/2}$ is uniformly continuous on $H$, by $(Hf)$.
\subsection{An extension of Theorem \ref{EM}}\label{emplus}
Our first result states that Theorem \ref{EM} can be enriched by
an additional uniformity in $h_n\le h\le \wth_n$ in the supremum
appearing in (\ref{ben}), provided that $(h_n)\suite$ and
$(\wth_n)\suite$ do fulfill assumption $(HV)$. We also refer to
Einmahl and Mason (2005), where the authors provided
some consistency results for kernel type function estimators that
hold uniformly in the bandwidth (see also Varron (2008) for an
improvement in the case of kernel density estimation).
\begin{theo}\label{T1}
Assume that $(HG)$, $(H f)$, $(H\CC)$ and $(HK1)-(HK3)$ are
satisfied. Let $(h_n)\suite$ and $(\wth_n)\suite$ be two sequences
of constants fulfilling $(HV)$ as well as $h_n=o(\wth_n)$. Then we have almost surely \beq
\limn \sup_{z\in H,\; g\in \GG,\;h_n\le h\le \wth_n} \frac{\mid
W_n(g,h,z)\mid}{\sqrt{2 n h^d \log(h^{-d})}}=\Delta(\GG)\mmi
K\mmi_{\lab,2}.\eeq
\end{theo}
The proof of Theorem \ref{T1} is provided in Section \ref{PT1}.
\begin{rem}
Einmahl and Mason (2005) have proved a result strong
enough to derive that, under weaker conditions than those of
Theorem \ref{T1}, we have almost surely \beq
\lsn\mathop{\sup_{z\in H,\;g\in \GG,}}_{h\in [\frac{c\log
n}{n},1]}
\frac{f_Z(z)^{1/2}W_n(g,h,z)}{\sqrt{nh^d\log(1/h)+\log\log
n}}<\infty.\label{rlp}\eeq However, the finite constant appearing
on the right hand side of (\ref{rlp}) is not explicit in their
result. The main contribution of Theorem \ref{T1} is that the
right hand side of (\ref{rlp}) is explicit, by paying the price of
making stronger assumptions.
\end{rem}
\begin{rem}
As Theorem \ref{T1} is an extension of Theorem \ref{EM} of Einmahl
and Mason, all the corollaries of Theorem \ref{EM} (see
Einmahl and Mason (2000)) can be enriched with a uniformity in the
bandwidth.
\end{rem}
\subsection{Some applications of Theorem \ref{T1} to data-driven bandwidth selection}
The main statistical interest of Theorem \ref{T1} is that we can
derive the limit behavior of kernel regression estimators with
data-driven bandwidth. Let us consider such a random bandwidth
$\ovh_n(z)=h(z,Y_1,\ldots,Y_n,Z_1,\ldots,Z_n)$ that depends on the
sample as well as on the point $z\in \RRR^d$. In the sequel, $Id$
shall denote the identity function. Our next corollary gives the
a.s. limit behavior of the Nadaraya-Watson estimator
$$ r_n(z) = \sum_{i=1}^n \frac{K\Big(\frac{z-Z_i}{h}\Big)}{\sum_{j=1}^n K\Big(\frac{z-Z_j}{h}\Big)} Y_i
$$
of the regression function $r(z):=\EEE(Y\mid Z=z)$, when
$\ovh_n(\cdot)$ satisfies some mild conditions. Note that the
asymptotics are given for $r_n(\cdot)-\ovr(\ovh_n,\cdot)$, with
$$\ovr(\ovh_n,z):= \ovh_n^{-1} \ili_{\RRR^d\times
\RRR^k}yK\poo\frac{u-z}{\ovh_n}\pff f_{Y,Z}(y,u)dudy.$$ The
\textit{random} differences $\ovr(\ovh_n,z)-r(z)$ can be
controlled by analytic arguments  as soon as the a.s. limit behavior is
known.
\begin{coro}\label{data}
Assume that $\ovh_n(\cdot)$ satisfies almost surely (resp. in
probability)
$$0<\lin \frac{\log(1/\ovh_n)}{\log n}\le \lsn \frac{\log(1/\ovh_n)}{\log
n}<1.$$ Then, we have
$$\lsn \sup_{z\in
H}\frac{\pm
\sqrt{f_Z(z)}(r_n(z)-\ovr(\ovh_n,z)) \sqrt{n\ovh_n(z)^d}}{\sqrt{2\Delta(Id,z)\log(\ovh_n(z)^{-d})}}=\mmi
K\mmi_{\lab,2},$$ almost surely (resp. in probability).
\end{coro}
\textbf{Proof}: The proof involves continuity arguments for
$\Delta(Id,\cdot)$ and the fact that the numerator and denominator
of $r_n(z)$ are specific forms of the general object $W_n$
appearing in Theorem \ref{T1}. We also consider the countable
collection of events
$$\aoo n^{-r}\le \ovh_n \le n^{-r'}\text{ for all large }
n\aff,\;r,r'\in  \QQQ\cap(0,1).$$ On each of these countable
events, the sequence $\ovh_n$ can be bounded from below and above
by sequences $h_n$ and $\mathfrak{h}_n$ fulfilling condition
$(HV)$. We omit technical details.$\Box$ \\
\lb
 \textbf{Example 1:
Tsybakov's plug-in selection rule}:\lb Tsybakov
(1987) considered a plug-in bandwidth selection rule when
$d=k=1$. In that case, he suggested that, for a given point $z\in
\RRR$, the bandwidth should be chosen of the form
$$\ovh_n(z):=\hat{\beta}_n(z)n^{-1/5},$$
where $\hat{\beta}_n(z)$ is a consistent estimate of the
theoretical quantity $\beta(z)$ that minimizes the asymptotic
square error of $r_n(z)$. Under the conditions stated in
Tsybakov (1987), since most of them being consequences of the
assumptions of Theorem \ref{T1}, the plug-in bandwidth satisfies
the assumptions of Corollary \ref{data}.\\
\lb
\textbf{Example 2:
cross validation}:
\lb We again consider the case $d=k=1$. An
important example is the bandwidth $\ovh_n$ that minimizes the
sample-based quantity
$$CV(h):=\frac{1}{n} \sliin [Y_i-r_{n,-i}(Z_i)]^2w(Z_i),\; h\in
[n^{-1+\dd},n^{-\dd}],$$ where $w$ is a weight function on $\RRR$
and $\dd>0$ is a fixed (small) value. We refer to Clark
(1975) and Priestley and Chao (1972) for
more details on that technique. By construction the random
sequence $\ovh_n$ satisfies the assumptions of Corollary
\ref{data}. Moreover, it is shown in Härdle \textit{et al.} (1988) that, under
mild conditions, we have $\ovh_n\sim C_0 n^{-1/5}$ in probability,
for a theoretical constant $C_0$.
\subsection{Asymptotic confidence bands by empirical likelihood}\label{elplus}
Empirical likelihood methods in statistical inference have been
introduced by Owen (2001). This nonparametric
technique has suscitated much interest for several practical
reasons, the most important one being that it directly provides
confidence intervals without requiring further approximation
methods, such as the estimation of dispersion parameters.
Moreover, empirical likelihood is a very versatile tool which can
be adapted in many different fields, for instance in estimation of
densities or conditional expectations by kernel smoothing methods.
The idea can be summarised as follows : consider an independent,
identically distributed sample $(Y_i,Z_i)_{1\le i\le n}$ taking
values in $\RRR^{d'}\times \RRR^d$. Given $h>0$, $z\in H$, a
function $g$ from $\RRR^{d'}$ to $\RRR^k$ and a (kernel) real
function $K$, define the following centring parameter, which plays
the role of a deterministic approximation of $\EEE\po g(Y)\mid
Z=z\pf$ : \beq m(g,h,z):= \frac{\EEE\pooo g(Y)K\poo
\frac{Z-z}{h}\pff\pfff}{\EEE\pooo K\poo
\frac{Z-z}{h}\pff\pfff}.\label{mghz}\eeq This quantity is the root
of the following equation in $\tht$: \beq \EEE\pooo K\poo
\frac{Z-z}{h}\pff\poo g(Y)-\tht\pff\pfff=0,\eeq which naturally
leads to the following formula for a confidence interval (around
$m(g,h,z)$) by empirical likelihood methods (for more details see,
e.g., Owen (2001), chapter 5) : \beq I_n(g,h,z,c):=\{\tht
\in \RRR,\; \RR_n(\tht,g,h,z)\geq c\},\label{Igzhalp}\eeq where
$c\in (0,1)$ is a given critical value that has to be chosen in
practice, and where \beq \RR_n(\tht,g,h,z):=\max\aoo \proliin
np_i, \sliin p_i K\poo \frac{Z_i-z}{h}\pff \po
g(Y_i)-\tht\pf=0,\;p_i\geq 0,\;\sliin
p_i=1\aff.\label{RRnthtgz}\eeq It is known (see, e.g.,
(2001), chapter 5) that, for fixed $z\in \RRR^d$ and fixed
$g$, we can expect \beq m(g,h,z)\in I_n(g,h,z,c)\label{app}\eeq to
hold with probability equal to $\PPP(\chi^2\le -2\log c)$,
ultimately as $\nif,\;h\rar 0,\;nh^d\rar\infty$ (see e.g., Owen,
chapter 5). A natural arising question is:
\begin{itemize} \item
Can we expect (\ref{app}) to hold uniformly in $z,g$ and $h$?
\item In that case, how much uniformity can we get?
\end{itemize}
Uniformity in $g$ and $z$ would allow to construct asymptotic
\textit{confidence bands} (instead of simple confidence
intervals), while a uniformity in $h$ would allow more flexibility
in the practical choice of that smoothing parameter. Our Theorem 3
provides a tool strong enough to give some positive answers to
these questions. We shall focus on the case where
$\GG=\{1_C,\;C\in\CC\}$ for a class of sets $\CC$. We will also
make an abuse of notation, by identifying $\CC$ and $\GG$, and
hence, we shall write $m(C,h,z)$ for $m\po 1_C,h,z\pf$ and so on.
Write the conditional variance of $1_C(Y)$ given $Z=z$ as follows
: \beq\sig^2(C,z):=\PPP\po Y\in C\mid Z=z\pf-\PPP^2\po Y\in C\mid
Z=z\pf,\;C\in \CC,\;z\in H\label{sigCz}.\eeq
The next theorem
shows that we can construct, by empirical likelihood methods
(recall (\ref{Igzhalp})), confidence bands around the centring
parameters $m(C,h,z)$ with lengths tending to zero at rate
$\sqrt{2\sig^2(C,z)\log(h^{-d})/nh^d}$ when $\nif$ and $h_n\le h
\le \wth_n$. We make the following assumptions on $h_n$, $\wth_n$
and $\CC$ :
\begin{tabbing}
$\hskip5pt$ \= $(HG')$$\hskip5pt$ \= $\CC$ is a VC class
satisfying
$\displaystyle{\inf_{z\in H}\inf_{C\in \CC} \sig^2(C,z)=:\beta>0}$.\\
 \>$(HV')$ \> $h_n \downarrow 0,\;nh^d_n
\uparrow\infty,\;\log(1/h_n)/\log_2 n\rar \infty,\;
nh^d_n/\log(1/\wth_n)\rar\infty.$
\end{tabbing}
Note that $(HV')$ is equivalent to $(HV)$ in the specific case
where $p=\infty$.
\begin{theo}\label{T2}
Under assumptions $(Hf)$, $(HK1)-(HK3)$,  $(HG')$ and $(HV')$, as
well as $h_n=o(\wth_n)$, we have almost surely:
\beq
\limn
\mathop{\sup_{z\in H,\;C\in \CC,}}_{h_n\le h\le \wth_n} \frac{-\log
\RR_n\poo m(C,h,z),C,h,z\pff}{\log(h^{-d})}=1.\label{conv}
\eeq
\end{theo}
The proof of Theorem \ref{T2} is provided in Section \ref{PT2}.
\begin{rem}
Theorem \ref{T2} implies that, for an arbitrary $\e>0$, taking
$c=h^{d+\e}$ when constructing confidence regions as in
(\ref{Igzhalp}) ensures that each $m(C,h,z)$ belongs to its
associated confidence interval $I_n(C,h,z,c)$. Moreover, this
claim turns out to be false when taking $c=h^{d-\e}$ with $\e>0$.
This shows that one cannot go below the theoretical limit $c=h^d$
without loosing uniformity in $C,h$ and $z$.
\end{rem}
\begin{rem}
In order to obtain a confidence band for $m(C,h,z)$ uniformly in
$C, h$ and $z$, we need the limiting distribution of
\begin{eqnarray} \label{supp}
\sup_{z,C,h} \Big[-\log \RR_n\poo m(C,h,z),C,h,z\pff\Big]/\log(h^{-d}),
\end{eqnarray}
so Theorem 3 is not
sufficient for this.   Obtaining such a limit law is a real challenge
in itself, and is beyond the scope of this paper.  We leave that
problem as an open problem.  In the case of univariate kernel density
estimation, Bickel and Rosenblatt (1973) showed that the supremum
over the transformed kernel density estimator, obtained after a proper
rescaling and a proper translation, converges to an extreme value
distribution.  The simulations in Section 2 suggest that a proper
linear transformation of $[-\log
\RR_n (m(C,h,z),C,h,z)]/\log(h^{-d})$ (depending on $z, C$ and
$h$) might also lead to a nondegenerate limiting distribution.
\end{rem}

\section{Simulation results}\label{simulations}
\begin{figure}
\centering
\includegraphics{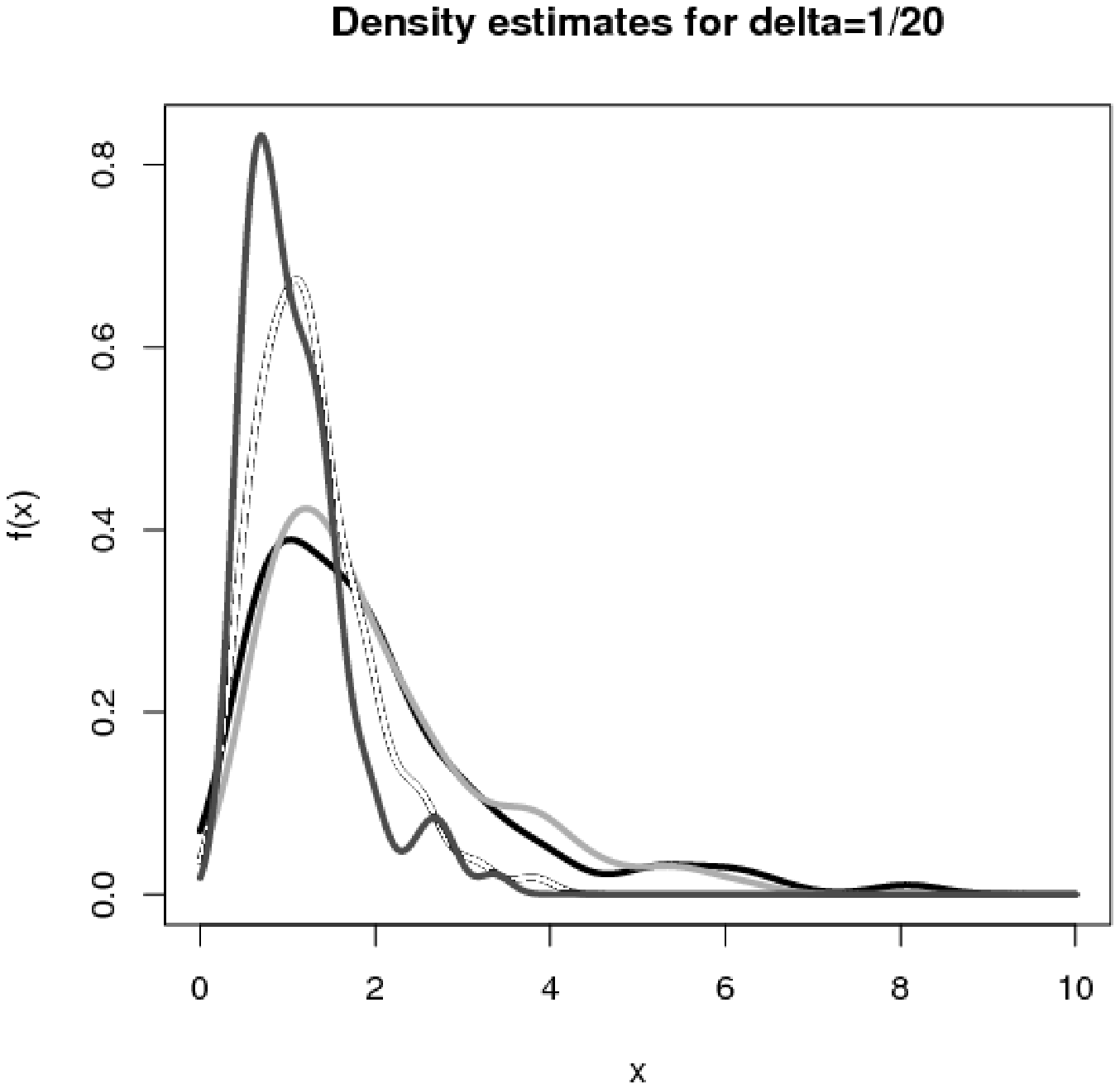}
\isucaption{Estimated densities of the supremum in (\ref{supp}) for $\delta=1/20$.
The black curve corresponds to $n=50$, the light gray curve to $n=100$, the white
curve to $n=500$ and the dark gray curve to $n=1000$.}
\label{fig1}
\end{figure}

\begin{figure}
\includegraphics{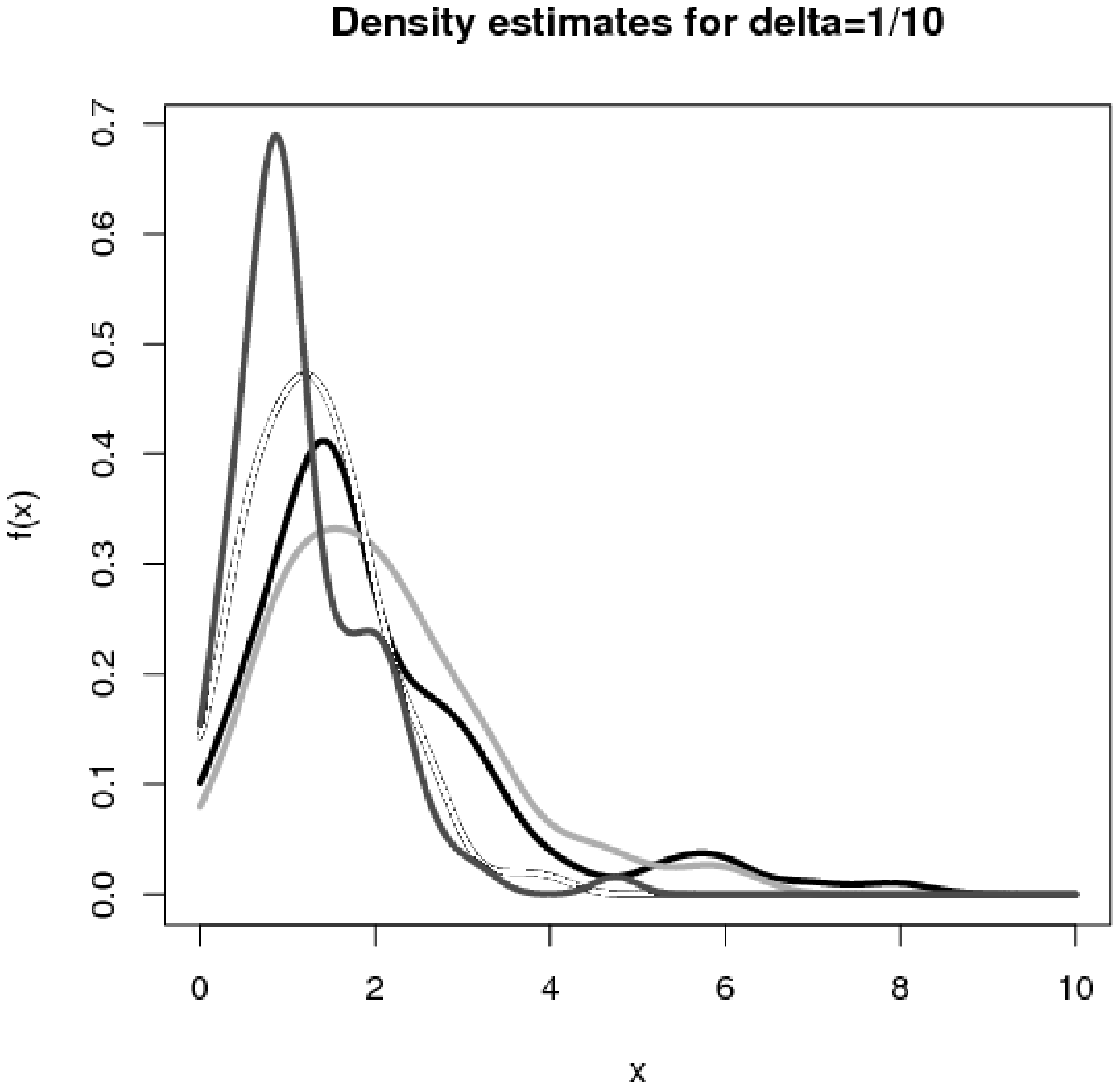}
\isucaption{Estimated densities of the supremum in (\ref{supp}) for $\delta=1/10$.
The light gray curve corresponds to $n=50$, the black curve to $n=100$, the white
curve to $n=500$ and the dark gray curve to $n=1000$.}
\label{fig2}
\end{figure}

A simulation study is carried out to illustrate the convergence
stated in (\ref{conv}). We estimate the density of (\ref{supp})
for four different sample sizes:
$n=50,100,500,1000$. We specified the following parameters:
\begin{enumerate}
  \item $Z$ is uniformly distributed on $[0,1]$. Given $Z=z$, $Y$
  has an exponential distribution with expectation $1/z$.
  \item $C$ is the class of intervals $[0,t],\;t\in [1,2]$.
  \item $H=[0.25,0.75]$.
  \item $h_n=n^{-1/5-\dd}$ and $\wth_n=n^{-1/5+\dd}$, with $\dd=1/20$.
\end{enumerate} For each sample size, the density is estimated as
follows :
\begin{itemize}
\item 100 independent samples are simulated (which is enough
since the density is univariate).
\item For each sample, the supremum in (\ref{supp})
is approximated by a maximum over a finite grid of
size 50.
\item Finally, the density of (\ref{supp}) is estimated by using a
Parzen-Rosenblatt density estimator, applied to the 100 obtained
values.  We used an Epanechnikov kernel and the bandwidth was obtained
from cross validation.
\end{itemize}
Figure \ref{fig1} shows the density estimates for
$n=50,100,500,1000$. Figure \ref{fig2} has been obtained
from a second simulation study, where the interval
$[h_n,\wth_n]$ has been widened ($\dd=1/10$).
As already mentioned in Remark 1.4, Figures \ref{fig1} and \ref{fig2}
suggest that after a proper linear transformation, the distribution of
(\ref{supp}) might converge to a non-degenerate limiting distribution.

\section{Proof of Theorem \ref{T1}}\label{PT1}
For ease of notations, we just prove Theorem \ref{T1} when $k=1$.
A close look at the proof shows that there is no loss of
generality assuming $k=1$.
\subsection{Truncation}
 We start our proof of Theorem \ref{T1} as
Einmahl and Mason did in their proof of Theorem \ref{EM}. As the
support of $K$ is bounded and as $\wth_n\rar 0$ we have almost
surely, for all large $n$ and for all $z\in H,\;g\in \GG,\;h_n\le
h\le \wth_n$, \begin{align} \nono
W_n(g,h,z)=&f_Z(z)^{-1/2}\sliin\cooo \poo
c_g(z)g(\wtY_i)+d_g(z)\pff
K\poo \frac{Z_i-z}{h}\pff\\
&\;-\EEE\pooo \poo c_g(z)g(\wtY_i)+d_g(z)\pff K\poo
\frac{Z_i-z}{h}\pff\pfff\cfff,\end{align} where
$\wtY_i:=Y_i1_{O'}(Z_i)$. Hence, we can suppose that
$Y_i=Y_i1_{O'}(Z_i)$ without changing the limiting behaviour of
the processes we are studying here. Now consider a sequence of
constants $(\gam_n)\suite$ fulfilling \beq \lin
\frac{\gam_n}{(n/\log(1/h_n))^{1/p}}>0,\label{abel}\eeq and
consider the truncated expressions, with $G$ denoting a measurable
envelope function of $\GG$ fulfilling $(HG)$,\begin{align} \nono
W_n^{\gam_n}(g,h,z):=&f_Z(z)^{-1/2}\sliin \cooo\poo
c_g(z)g(Y_i)1_{\{G(Y_i)\le \gam_n\}}+d_g(z)\pff K\poo
\frac{Z_i-z}{h}\pff\\
&\;-\EEE\pooo\poo c_g(z)g(Y_i)1_{\{G(Y_i)\le \gam_n\}}+d_g(z)\pff
K\poo \frac{Z_i-z}{h}\pff\pfff\cfff.\end{align} The following
lemma allows us to study these truncated versions of the
$W_n(g,h,z)$.
\begin{lem}\label{truncation}
Under the assumptions of Theorem \ref{T1} and under (\ref{abel})
we have almost surely:\beq \limn \sup_{g\in \GG,\;z\in H,\;h_n\le
h\le \wth_n}\frac{\Mid
W_n^{\gam_n}(g,h,z)-W_n(g,h,z)\Mid}{\sqrt{2nh^d\log(h^{-d})}}=0.\eeq
\end{lem}
\textbf{Proof}: A careful reading of the proof of Lemma 1 in
Einmahl and Mason (2000) shows that their assertions
(2.8) and (2.9) remain true after adding a uniformity in $g\in
\GG$ and $h_n\le h\le \wth_n$, which readily implies Lemma
\ref{truncation}. Note also that Lemma \ref{truncation} is obvious
when $(HG)$ is fulfilled with $p=\infty$. $\Box$\lb The two next
subsections are devoted to proving respectively the outer and
inner bounds of Theorem \ref{T1}.
\subsection{Outer bounds}\label{outer}
Fix $\e>0$. Our goal in this subsection is to show that, almost
surely \beq \lsn \sup_{g\in \GG,\;z\in H,\;h_n\le h\le
\wth_n}\frac{\mid W_n(g,h,z)\mid}{\sqrt{2nh^d\log(h^{-d})}}\le
\Delta(\GG)\mmi K\mmi_{\lab,2}(1+4\e).\label{obj1}\eeq To this
aim, we shall first discretise each of the sets $H$,
$[h_n,\wth_n]$ and $\GG$ into properly chosen finite grids, then
we shall control the oscillations between elements of the grids by
a combination of a concentration inequality which is due to
Talagrand (see also Massart (1989), Bousquet (2002) and Klein (2002) for sharpened
versions) and of an upper bound for the first moment of these
oscillations which is due to Einmahl and Mason (2000).
\subsubsection{Step 1: discretisations}
Consider three parameters $\dd_1\in (0,1)$, $\dd_2\in (0,1)$ and
$\rho\in (1,2)$ that will be chosen small enough in the sequel,
and define the following subsequence\beq \nk:= \coo \exp \po
k/\log k\pf\cff,\;k\geq
5,\;\;\;\;\;\;\;N_k:=\{\nkm,\nkm+1,\ldots,\nk-1\}.\label{nk}\eeq
Note that $\nk/\nkm\rar 1$ and\beq \log\log \nk=\log
k(1+o(1)),\;\kif.\label{expo}\eeq We then construct the following
finite grid for each $k\geq 1$
\begin{align}
h_{\nk,R_k}:=\wth_{\nkm},\;\; \hkl:=&\rho^l \hk,\;l=0,\ldots,
R_k-1,\label{discret1}
\end{align}
where $R_k:= [\log(\wth_{\nkm}/\hk)/\log(\rho)]+1$, and $[u]$
denotes the only integer $q$ fulfilling $q\le u<q+1$. Denote by
$\mid z\mid_d:= \max_{i=1,\ldots,d}\mid z_i \mid$ the usual
maximum norm on $\RRR^d$. Now, for fixed $k$ and $0\le l\le R_k$,
we construct a finite grid $\MM_{k,l}\subset H$ such that, given
$z\in H$, there exists $\zklj\in \MM_{k,l}$ fulfilling $\mid z-
\zklj \mid_d<\dd_1\hkl$. Note that one can construct this grid so
as $\sharp\MM_{k,l}\le C(\dd_1\hkl)^{-d}$, where $C$ is a constant
that depends only on the volume of $H$. Now set $\gam_n:=\dd_2 \po
\nk/\log(1/\hk^d)\pf^{1/p}$, for each $k\geq 5$, $n\in N_k$. By
Lemma \ref{truncation}, showing (\ref{obj1}) is equivalent to
showing that \beq \lsn \sup_{z\in H,\;g\in\GG,\;h_n\le h\le
\wth_n}\frac{\mid
W_n^{\gam_n}(g,h,z)\mid}{\sqrt{2nh^d\log(h^{-d})}}\le
\Delta(\GG)\mmi K\mmi_{\lab,2}(1+4\e)\label{obj2}\eeq almost
surely, for a proper choice of $\dd_2>0$.
\subsubsection{Step 2: a
discrete version of (\ref{obj1})} Given a real function $\psi$
defined on a set $S$, we shall write: \beq
\mmi\psi\mmi_{S}:=\sup_{s\in S}\mid \psi(s)\mid.\eeq Recall that,
since $f_Z$ is bounded away from $0$ on $H$, we can define \beq
\gam:=\inf_{z\in H} f_Z(z)>0.\label{gam}\eeq Also write, for
convenience of notations \beq \mmi c\mmi_{\GG\times H}:=
\sup_{g\in \GG,\;z\in H} \mid c_g(z)\mid,\;\;\; \mmi
d\mmi_{\GG\times H}:= \sup_{g\in \GG,\;z\in H} \mid
d_g(z)\mid.\eeq
 Our first lemma is a version of (\ref{obj2}) which is discretised along the finite grids defined in Step 1.
\begin{lem}\label{dis}
For any choice of \beq 0<\dd_2< \e\gam^{1/2}\Delta(\GG)\mmi
K\mmi_{\lab,2}/(6\mmi c\mmi_{H\times \GG}\mmi
K\mmi_{\RRR^d}),\label{rouf}\eeq for any finite collection
$\{g_1,\ldots,g_q\}\subset \GG$ and for any $\dd_1\in
(0,1),\;\rho\in (1,2)$, we have \beq \lsk \mathop{\max_{n\in
N_k,\;1\le \ell \le q,}}_{0\le l\le R_k,\;\zklj\in
\MM_{k,l}}\frac{\mid
W_n^{\gam_n}(g_\ell,\hkl,\zklj)\mid}{\sqrt{2\nk\hkl^d\log(1/\hkl^d)}}\le
\Delta(\GG)\mmi K\mmi_{\lab,2}(1+\e).\eeq
\end{lem}
\textbf{Proof}: We can assume here that $q=1$ with no loss of
generality. We rename in this proof $g_1$ to $g$. We define, for
$z\in H,\;h>0$ and $\ovg \in \GG$, \beq \psi_{\nk,h,z,\ovg}:\;
(y,x)\mapsto f_Z(z)^{-1/2}\coo c_{\ovg}(z)\ovg(y)1_{\{G(y)\le
\gam_{\nk}\}}+ d_{\ovg}(z)\cff K\poo\frac{x-z}{h}\pff.\eeq First
note that, for each $k\geq 5,\;0\le l\le R_k$ and $\zklj \in
\MM_{k,l}$, we have
\begin{align} \nono \mmi \psi_{\nk,\hkl,\zklj,g}\mmi_{\RRR^{d'}\times\RRR^d} \le & \po
\mmi c\mmi_{H\times \GG} \gam_{n_k}+\mmi d \mmi_{H\times \GG}\pf \gam^{-1/2}\mmi K \mmi_{\RRR^d}\\
\le &2\mmi c\mmi_{H\times \GG} \gam^{-1/2}\mmi K \mmi_{\RRR^d}
\dd_2
(\nk\hk^d/\log(1/\hk^d))^{1/2}\label{pasta}\\
\le & \frac{\e}{3}\mmi
K\mmi_{\lab,2}^2\Delta(\GG)(\nk\hk^d/\log(1/\hk^d))^{1/2}\label{paste},
\end{align}
where (\ref{pasta}) holds for all large $k$, uniformly in $0\le
l\le R_k$ and $\zklj\in \MM_{k,l}$, according to assumption
$(HV)$, and where (\ref{paste}) holds by (\ref{rouf}). Moreover
we have (recall $(HK2)$)
\begin{align}
\nono&\Var\poo \psi_{\nk,\hkl,\zklj,g}(Y,Z)\pff\\
\nono \le &\;\EEE\poo
\psi^2_{\nk,\hkl,\zklj,g}(Y,Z)\pff\\
\nono \le & \;\EEE\pooo f_Z(\zklj)^{-1}\poo
c_g(\zklj)g(Y)+d_g(\zklj)\pff^2 K\poo
\frac{Z-\zklj}{\hkl}\pff^2\pfff\\ &\;+ f_Z(\zklj)^{-1}\mmi
d\mmi_{H\times \GG}^2 \mmi K\mmi_{\RRR^d}^2 \PPP \poo \ao G(Y)\geq
\gam_{\nk}\af\cap \ao\mid Z-\zklj\mid_d\le
\hkl/2\af\pff\label{vol}\\
\nono =:& A_{1,\zklj}+A_{2,\zklj}.
\end{align}
The first term on the right hand side of (\ref{vol}) is equal to
\begin{align}
\nono A_{1,\zklj}=&\ili_{\mid z-\zklj\mid_d\le \hkl/2}\EEE\pooo
\po c_g(\zklj)g(Y)+d_g(\zklj)\pf^2 \Mid Z=z\pfff
\frac{f_Z(z)}{f_Z(\zklj)} K^2\poo\frac{z-\zklj}{\hkl}\pff dz.
\end{align}
It follows, by making use of assumption $(H\CC)$, that there
exists a function $r(\cdot)$ fulfilling $r(u)\rar 0$ as $u\rar 0$
and such that
\begin{align} A_{1,\zklj}\le & \ili_{\mid z-\zklj\mid_d\le \hkl/2}\Delta^2(g,z)
\frac{f_Z(z)}{f_Z(\zklj)}
K^2\poo\frac{z-\zklj}{\hkl}\pff dz+r(\hkl)\\
\le &\Delta^2(\GG)\hkl^d\ili_{[-1/2,1/2]^d}K^2(u)\frac{f_Z\po
\zklj+\hkl u\pf}{f_Z(\zklj)} du\dot\po1+r(\hkl)\pf\\
\le  &\Delta^2(\GG)\mmi K\mmi_{\lab,2}^2 \hkl^d(1+
\varepsilon_{k,l}),
\end{align}
where \beq \varepsilon_{k,l}:=\sup_{z\in H,\;\mid u \mid_d\le
1/2}\Bigg |\frac{f_Z\po z+\hkl u\pf}{f_Z(z)}\po1+r(\hkl)\pf-1
\Bigg |.\eeq By assumption $(H f)$ and since $\hkl\le
\wth_{\nkm}\rar 0$ we readily infer that
$$\limk \max_{0\le l\le R_k}\varepsilon_{k,l}= 0.$$ Moreover we have, uniformly in $0\le l\le R_k$ and $\zklj\in \MM_{k,l}$
(recall $(HG)$ and $(Hf)$)
\begin{align}
\nono&\PPP \poo \{G(Y)\geq \gam_{\nk}\}\cap \{\mid
Z-\zklj\mid_d\le\hkl/2\}\pff\\
\nono\le &\gam_{\nk}^{-2}\ili_{\mid z-\zklj\mid_d\le
\hkl/2}\EEE\poo G^2(Y)\mid Z=z\pff f_Z(z)dz\\
\nono\le & \gam_{\nk}^{-2}
\hkl^d\;\alp^2\ili_{[-1/2,1/2]^d}f_Z\po\zklj+\hkl u\pf du.\\
\nono\le & \gam_{\nk}^{-2} \hkl^d\;\alp^2\mmi f_Z\mmi_{O}.
\end{align}
 As $\gam_{\nk}\rar \infty$ we conclude that,
for all large enough $k$ and for each $0\le l\le R_k$, $\zklj\in
\MM_{k,l}$,
 \beq \Var\poo \psi_{\nk,\hkl,\zklj,g}(Y,Z)\pff\le
\Delta^2(\GG)\mmi K\mmi_{\lab,2}^2(1+\e)\hkl^d.\label{var1}\eeq
Given a real function $g: \RRR\times\RRR^d\mapsto \RRR$, we shall
write \beq T_n(g):=\sliin \aooo g(Y_i,Z_i)-\EEE\poo
g(Y_i,Z_i)\pff\afff.\label{Tn}\eeq Combining (\ref{paste}) and
(\ref{var1}) making use of the maximal version of Bernstein's
inequality (see, e.g. Einmahl and Mason (1996), Lemma
2.2) repeatedly for each $0\le l \le R_k,\;\zklj\in \MM_{k,l}$, we
have, for all large $k$ (recall that $\sharp \MM_{k,l}\le C
\dd_1^{-d} \hkl^{-d}$),
\begin{align}
\nono&\PPP\pooo\mathop{\max_{n\in N_k,\;0\le l\le R_k,}}_{\zklj\in
\MM_{k,l}}\frac{\mid
W_n^{\gam_n}(g,\hkl,\zklj)\mid}{\sqrt{2\nk\hkl^d\log(1/\hkl^d)}}>
\Delta(\GG)\mmi K\mmi_{\lab,2}(1+\e)\pfff\\
\nono \le& \sli_{l=0}^{R_k} \sharp \MM_{k,l} \max_{\zklj\in
\MM_{k,l}} \PPP\pooo \max_{n\in N_k}
\Mid T_n\po \psi_{n_k,h_{n_k,l},\zklj,g}\pf \Mid\geq \Delta(\GG)\mmi K\mmi_{\lab,2}(1+\e)\sqrt{2\nk\hkl^d\log(1/\hkl^d)}\pfff\\
\nono\le& \sli_{l=0}^{R_k} \frac{C}{\dd_1^d\hkl^d} 2\exp \poo- \po
1+\e^2/(1+\e)\pf \log (1/\hkl^d)\pff\\
\nono\le & \frac{2C}{\dd_1^d}\sli_{l=0}^{R_k} \hkl^{d\e^2/2} \\
=& \frac{2C}{\dd_1^d} \sli_{l=0}^{R_k}
\rho^{ld\e^2/2}\hk^{d\e^2/2} = \frac{2C}{\dd_1^d}\hk^{d\e^2/2}
\frac{\rho^{(R_k+1)d\e^2/2}-1}{\rho^{d\e^2/2}-1} \le
\frac{2C\rho^{d\e^2/2}}{\dd_1^d(\rho^{d\e^2/2}-1)}\wth_{\nkm}^{d\e^2/2},\label{rtyq}
\end{align}
where the last inequality is a consequence of $R_k:=
[\log(\wth_{\nkm}/\hk)/\log(\rho)]+1$. As
$\log(1/\wth_{\nkm})/\log\log \nkm\rar \infty$ (assumption (HV)),
and by (\ref{expo}), the right hand side of expression
(\ref{rtyq}) is summable in $k$. The proof of Lemma \ref{dis} now
readily follows by making use of the Borel-Cantelli lemma.
$\Box$\lb
\subsubsection{Step 3: end of the proof of Theorem \ref{T1}} Our
next lemma allows us to extend the uniformity in Lemma \ref{dis}
to the whole sets $\GG$, $[\hk,\wth_{\nkm}]$ and $H$, provided
that $\dd_1>0,\dd_2>0,\rho>1$ and $\{g_1,\ldots,g_{q}\}$ have been
properly chosen. Before stating our lemma, we need to recall three
facts. We shall be able to properly discretise the class $\GG$ by
making use of the following result, which is a straightforward
adaptation of Lemma 6 of Einmahl and Mason (2000).
\begin{fact}[Einmahl, Mason, 2000]\label{osso}
Given $\varepsilon>0$, there exists $h_{0,\varepsilon}>0$ and a
finite subclass $\{g_1,\ldots,g_q\}\subset \GG$ (that may depend
on $\varepsilon$) fulfilling \beq
\mathop{\sup_{0<h<h_{0,\varepsilon},}}_{z\in H,\;g\in
\GG}\;\min_{\ell=1,\ldots,q} h^{-d}f_Z(z)^{-1}\EEE \cooo \poo \po
c_g(z)g(Y)+d_g(z)\pf-\po c_{g_\ell}(z)g_\ell(Y)+
d_{g_\ell}(z)\pf\pff^2 K^2\poo\frac{Z-z}{h}\pff \cfff\le
\varepsilon/2.\nono\eeq
\end{fact}\vk
Now define the following distances on $\GG$: \begin{align}\nono
d^2(g_1,g_2):=& \mathop{\sup_{0<h<h_{0,\varepsilon},}}_{z\in H}
h^{-d}f_Z(z)^{-1}\EEE \coooo \poo \po
c_{g_1}(z)g_1(Y)+d_{g_1}(z)\pf-\po c_{g_2}(z)g_2(Y)+
d_{g_2}(z)\pf\pff^2\\
\nono&\times K^2\poo\frac{Z-z}{h}\pff \cffff
,\\
\wtd(g_1,g_2):=& \max\ao d(g_1,g_2),\;\mmi
c_{g_1}-c_{g_2}\mmi_H,\;\mmi d_{g_1}-d_{g_2}\mmi_H\af\label{dist1}
.\end{align} We write $\mid K\mid_v$ for the total variation of
$K$ and we set, for $\psi:\RRR^d\mapsto \RRR$,
\begin{align}\mathrm{\omega}_{\psi}(\dd):=& \sup_{z_1,z_2\in
H,\;\mid z_1-z_2 \mid_d \le \dd} \Mid \frac{\psi (z_2)}{f_Z(z_2)}-\frac{\psi(z_1)}{f_Z(z_1)}\Mid,\;\dd>0,\\
\beta_1:=&\sup_{z\in O} \EEE\poo \po G^2(Y)+1\pf\Mid
Z=z\pff<\infty,\\
B:=& 4\beta_1 \mmi f_Z\mmi_{O}\mmi f_Z^{-1}\mmi_O\poo \mmi
K\mmi^2_{\RRR^d}+\po\sup_{g\in \GG}\mmi c_g
\mmi_{O}^2\vee\sup_{g\in \GG}\mmi d_g \mmi_{O}^2\pf \mid
K\mid_v^2\pff.
\end{align}
The following fact is a straightforward adaptation of Lemma 4 and
Lemma 6 in (2000).
\begin{fact}[Einmahl, Mason, 2000]\label{f1}
Fix $\varepsilon>0$. For any $\dd\in (0,1/2)$ and
$0<h<h_{0,\varepsilon}$ fulfilling \beq z+ (2h) u\in O\text{ for
each }z \in H\text{ and for each }u \in \RRR^d\text{ with }\mid u
\mid_d \le 1 \label{arg}\eeq and for all large $k$ we have, for
each $\rho\in (1,2]$, $z_1,z_2\in H$ with $\mid z_1-z_2\mid_d \le
(\dd h)$, and for each $g_1,g_2\in \GG$ fulfilling
$\wtd^2(g_1,g_2)\le \varepsilon$,
\begin{align}\nono &\EEE\pooo \poo\psi_{\nk,z_1,\rho
h,g_1}(Y,Z)-\psi_{\nk,z_2,h,g_2}(Y,Z)\pff^2\pfff\\
\le& B \poo \mathrm{\omega}_{c_{g_2}}^2(\dd h)\vee
\mathrm{\omega}_{d_{g_2}}^2(\dd h)+\rho-1+\dd+\varepsilon \pff
h^d.\label{qalp}\end{align}
\end{fact}
\textbf{Remarks}: Assumption $(\ref{arg})$ is just technical, in
order to have the continuity arguments of Einmahl and Mason valid.
The presence of the term $\rho-1$ on the right hand side of
(\ref{qalp}) is due to the fact that we take care of the
differences $h/\hkl-1$, which are implicitly handled in Lemma 6 of
Einmahl and Mason (2000).\lb The third fact is also
largely inspired by the ideas of Einmahl and Mason
(2000). We remind that the uniform entropy number of a
class of functions $\FF$ with measurable envelope $F$ is defined
as \beq\nono \NN(\e,\FF):=\sup_{Q\text{ proba}} \min \ao p\geq
1,\;\exists (g_1,\ldots,g_p)\in \FF^p,\; \sup_{g\in
\FF}\min_{i=1,\ldots,p} \mmi g-g_i\mmi_{Q,2}\le \e\mmi F\mmi_{Q,2}
\af,\eeq where the supremum is taken over all probability measures
$Q$. The following fact is proved in Varron (2008, Proposition
2.1).
\begin{fact}[Varron, 2008] \label{f2}
Let  \(\mcal{\FF}\) be a  class of functions on $\RRR^d$ with
measurable envelope function $F$ satisfying, for some constants
$\tau>0$ and $h\in (0,1)$,
\[\sup_{g\in\mcal{\FF}}\;\mathrm{Var}\po g(Z_1)\pf\le \tau^2 h^d
.\] Assume that there exists $\dd_0,C,v,\beta_0>0$ and $p>2$
fulfilling, for all $0<\e<1$,
\begin{align}
\mcal{N}(\e,\FF)\le& C\e^{-v},\\
\EEE\poo F(Y)^2\pff\le& \beta_0^2,\\
\sup_{g\in \FF,\;z\in \RRR^d} \mid g(z) \mid \le & \dd_0
(nh^d/\log(h^{-d}))^{1/p}.\label{hash3}
\end{align}
Then there exists a universal constant $A>0$ and a parameter
$D(v)>0$ depending only on $v$ such that, for fixed $\rho_0>0$, if
$h>0$ satisfies,
\begin{align}
&K_1:=\max\aoo 1,\po4\dd_0\sqrt{v+1}/\tau\pf^{\frac{1}{1/2-1/p}},\po\rho_0 \dd_0/\tau^2\pf^{\frac{1}{1/2-1/p}}\aff \le  \frac{nh^d}{\log(h^{-d})},\label{hash1}\\
&K_2:=\min \ao1/(\tau^2\beta_0),\tau^2\af\geq h^d, \label{hash2}
\end{align}
then we have
\[\PPP\poo\max_{1\le m \le n}\mmi T_m\mmi_{\mcal{F}} \geq (\tau+\rho_0)D(nh^d\log(h^{-d}))^{1/2}\pff\le 4\exp\Big(-A(\frac{\rho_0}{\tau})^2\log(h^{-d})\pff.\]
\end{fact}
We can now state our second lemma, which will conclude the proof
of the outer bounds of Theorem  \ref{T1}. Recall that $\e>0$ was
fixed at the very beginning of our proof (see Section
\ref{outer}).
\begin{lem}\label{osci}
There exists a finite class $g_1,\ldots,g_{q}\in \GG$ as well as
two constants $\rho_\e>1$ and $\dd_{1,\e}>0$ small enough such
that, for each $1< \rho \le \rho_\e$ and each
$0<\dd\le\dd_{1,\e}$, we have almost surely :
\begin{align} \nono \lsk& \mathop{\max_{n\in N_k,}}_{0\le l\le
R_k-1}\;\sup_{g\in \GG}\inf_{1\le \ell\le
q}\;\mathop{\sup_{z_1,z_2\in H,\;\mid z_1-z_2\mid <\dd,}}_{
\hkl\le h \le \rho \hkl}\frac{\Mid
W_n^{\gam_n}(g,\hkl,z_1)-W_n^{\gam_n}(g_{\ell},h,z_2)\Mid}{\sqrt{2
\nk \hkl^d \log(1/\hkl^d)}} \\
\le & \Delta(\GG)\mmi K\mmi_{\lab,2}\e.\end{align}
\end{lem}
\textbf{Proof :}\lb Consider the class
\begin{align}
\nono\GG':=&\aoo (y,z)\mapsto u_1\po c_{g_1}(z_1)
g_1(y)1_{\{G(y)\le
t\}}+d_{g_1}(z_1)\pf K\poo \frac{z-z_1}{h}\pff\\
\nono&\;-u_2\po c_{g_2}(z_2) g_2(y)1_{\{G(y)\le
t\}}+d_{g_2}(z_2)\pf K\poo
\frac{z-z_2}{\mathfrak{h}}\pff,\;z_1,z_2\in
\RRR^d,\;g_1,g_2\in \GG,\\
\nono &\;t\geq 0,\;(h,\mathfrak{h})\in (0,1)^2,\;u_1,u_2\in
\coo\inf_H f_Z^{-1/2},\sup_Hf_Z^{-1/2}\cff\aff.\end{align} Recall
that $\gam=\inf_{H}f$ and note that $\GG'$ admits the following
function as an envelope function: \beq G':\;(y,z)\mapsto
2\gam^{-1/2}\po \mmi c\mmi_{H\times\GG} G(y)+ \mmi d\mmi_{H\times
\GG} \pf \mmi K \mmi_{\RRR^d}. \label{env} \eeq Set
$\beta_4^2:=\EEE\po {G'}^2(Y,Z)\pf<\infty$ (the finiteness of
$\beta_4$ follows from $(Hf)$ and $(HG)$). By an argument very
similar to that used in Lemma 5 of Einmahl and Mason (2000) we readily infer that there exist $C>0$ and
$v>0$ fulfilling \beq \NN(\e,\GG')\le C\e^{-v},\;\e\in
(0,1]\label{cover}.\eeq Recalling the notations of Fact \ref{f2},
we set $\varepsilon=D(v)^{-1}(1+\sqrt{2/A})^{-1}\e\Delta(\GG)\mmi
K\mmi_{\lab,2}.$  By Fact \ref{osso} and by $(H \mcal{C})$, for
any $\varepsilon>0$, we can choose a finite subclass
$\{g_1,\ldots,g_q\}\subset \GG$ such that $\GG$ is included in the
finite reunion of the corresponding balls with $\wtd$-radius
smaller than $\varepsilon/2$. For fixed $k\geq 5,\;0\le l\le
R_k-1,\;1\le \ell \le q$ and $\dd>0$, define the following class
of functions:
\begin{align} \nono\GG_{k,l,q,\dd}:=&\aoo
\psi_{\nk,z_1,h,g}-\psi_{\nk,z_2,\hkl,g_\ell},\;z_1,z_2\in
H,\;\mid
z_1-z_2\mid \le \dd,\\
\nono&\; \wtd(g,g_{\ell})\le \varepsilon/2,\; \hkl\le h\le \rho
\hkl\aff\label{GGkldd}.
\end{align} Obviously we always have $\GG_{k,l,\ell,\dd}\subset \GG'$.
By inclusion, all the classes $\GG_{k,l,\ell,\dd}$ inherit
properties (\ref{env}) and (\ref{cover}). Moreover, proving Lemma
\ref{osci} is equivalent to showing that, almost surely \beq \lsk
\mathop{\max_{0\le l\le R_k-1,}}_{1\le \ell \le q}
\frac{\max_{n\in N_k}\mmi T_n
\mmi_{\GG_{k,l,\ell,\dd}}}{\sqrt{2\nk\hkl^d\log(1/\hkl^d)}}\le
\Delta(\GG)\mmi K\mmi_{\lab,2}\e\label{rfa2}.\eeq As $\hkl \le
\wth_{\nkm}\rar 0$ and by Fact \ref{f1}, we can choose
$\dd_{1,\varepsilon}>0$ and $\rho_\varepsilon\in (1,2)$ such that,
for each $\dd_1\in (0,\dd_{1,\varepsilon}),\rho\in
(1,\rho_\varepsilon)$, for all large $k$ and for all $0\le l\le
R_k-1$, \beq
\sup_{\psi\in\GG_{k,l,\ell,\dd_{1,\varepsilon}}}\hkl^{-d}\EEE\po\psi^2(Y,Z)\pf\le
\varepsilon^2.\eeq Recalling that $\hk\le \hkl \le \wth_{\nkm}$
and assumption $(HV)$, we can choose $k$ large enough so that each
class $\GG_{k,l,\ell,\dd_1}$ fulfills conditions (\ref{hash1}) and
(\ref{hash2}) with $\beta_0:= \beta_4$, $h:=\hkl$, $n:= \nk$,
$\tau:=\varepsilon$, $\rho:=\sqrt{2/A}\tau$ and $C,v$ appearing in
(\ref{cover}). Hence, we have, uniformly in $0\le l \le R_k-1$ and
$1\le \ell \le q$,
\begin{align}
\nono&\PPP\pooo \max_{n\in N_k} \mmi T_n
\mmi_{\GG_{k,l,\ell,\dd_{1,\varepsilon}}}>\Delta(\GG)\mmi
K\mmi_{\lab,2}\e\sqrt{2\nk\hkl^d\log(1/\hkl^d)}
\pfff\\
\nono\le & \PPP\poo \max_{n\in N_k} \mmi T_n
\mmi_{\GG_{k,l,\ell,\dd_{1,\varepsilon}}}\geq
D(v)(\tau+\rho)\sqrt{2n\hkl^d\log(1/\hkl^d)}\pfff\\
\nono\le& 4 \exp \poo -2\log(1/\hkl^d)\pff.\end{align} Now, by
Bonferroni's inequality we have, for all large $k$,
\begin{align}
\nono &\PPP\poooo \bculi_{l=0}^{R_k-1}\bculi_{j=1}^{J_l}\max_{n\in
N_k} \mmi T_n \mmi_{\GG_{k,l,\dd_{1,\varepsilon}}}>\Delta(\GG)\mmi
K\mmi_{\lab,2}\e\sqrt{2\nk\hkl^d\log(1/\hkl^d)}
\pffff\\
\nono \le & \sli_{l=0}^{R_k-1} 4\sharp \MM_{k,l} \hkl^{2d} \nono
\le \frac{4 C}{\dd_{1,\e}^d} \sli_{l=0}^{R_k-1}\hkl^d \nono\le
\frac{4C\dd_{1,\e}^{-d}}{\rho^d-1} \rho^{dR_k} \nono \le
\frac{4C\rho^d\dd_{1,\e}^{-d}}{\rho^d-1} \wth_{\nkm}.
\end{align}
As $\log(1/\hk)/\log\log(\nk)\rar \infty$ by (HV) and
(\ref{expo}), the proof of Lemma \ref{osci} is concluded by a
straightforward use of the Borel-Cantelli lemma. $\Box$\lb
Combining Lemmas \ref{dis} and \ref{osci} we get, for any choice
of $\dd_1,\;\dd_2>0$ and $\rho>1$ small enough, \beq \lsk
\mathop{\max_{n\in N_k,}}_{0\le l\le R_k-1}\;\mathop{\sup_{z\in
H,}}_{\hkl\le h \le \rho\hkl}\frac{\Mid
W_n^{\gam_n}(g,h,z)\Mid}{\sqrt{2 \nk \hkl^d \log(1/\hkl^d)}}\le
\mmi K\mmi_{\lab,2}\Delta(\GG)(1+3\e)\;\text{a.s.}\label{obj3}\eeq
Now assertion (\ref{obj1}) is almost proved, provided that we
substitute $n_k \hkl^d \log(1/\hkl^d)$ by $nh^d\log(1/h^d)$ in the
LHS of (\ref{obj3}) at the minor cost of replacing $1+3\e$ by
$1+4\e$ in the RHS of (\ref{obj3}). This can be achieved by
noticing that $\max\{\mid \sqrt{n/\nk}-1\mid,\;\nkm<n\le\nk\}\rar
0$ together with the following assertion \beq
\lsk\mathop{\sup_{1<\rho'\le \rho,}}_{h\in (\hk,\wth_{\nkm})} \Mid
\sqrt{\frac{(\rho h)^d\log\po(\rho h)^{-d}\pf}{h^d\log(h^{-d})}}-1
\Mid \le \e/(1+5\e),\label{softmoutonorg}\eeq which, by routine
computations, turns out to be true if we choose $\rho>1$ small
enough. This concludes the proof of the outer bounds of Theorem
\ref{T1}.$\Box$
\subsection{Inner bounds}
Proving the inner bounds of Theorem \ref{T1} is a simple
consequence of Theorem \ref{EM}, since, almost surely,
\begin{align} \lin \mathop{\sup_{z\in H,\;C\in \CC,}}_{h_n\le h\le \wth_n}
\frac{\mid W_n(g,h,z)\mid}{\sqrt{2nh^d\log(h^{-d})}} \geq & \lin
\mathop{\sup_{z\in H,}}_{C\in \CC} \frac{\mid
W_n(g,h_n,z)\mid}{\sqrt{2nh^d_n\log(h^{-d}_n)}}= \Delta(\GG)\mmi
K\mmi_{\lab,2},\label{hok}
\end{align}
where (\ref{hok}) is a consequence of Theorem \ref{EM}.
\section{Proof of Theorem \ref{T2}}\label{PT2}
Our proof of Theorem \ref{T2} is inspired by chapter 5 in Owen (2001) and borrows some ideas of Chen \textit{et
al.} (2003). Set, for $n\geq 1$,
$C\in \CC$, $h>0$ and $z\in H$, \begin{align} X_n(C,h,z):= &\sliin
K\poo \frac{Z_i-z}{h}\pff \poo 1_C(Y_i)- m(C,h,z)\pff,\label{Xngz}\\
 S_n(C,h,z):=&f_Z(z)^{-1}
 \sliin {\cooo K\poo\frac{Z_i-z}{h}\pff \poo
1_C(Y_i)-m(C,h,z)\pff\cfff}^2,\label{Sngz}\\
w_{i,n}(C,h,z):=&\; K\poo \frac{Z_i-z}{h}\pff \poo 1_C(Y_i)-
m(C,h,z)\pff\label{wigz}.
\end{align}
The proof of Theorem \ref{T2} consists in showing that the
quantities
$$-2\log\poo\RR_n\po m(C,h,z),C,h,z\pf\pff,\;C\in \CC,\;z\in
H,\;h\in [h_n,\wth_n]$$ are asymptotically equivalent to \beq
U_n(C,h,z):=\frac{X_n(C,h,z)^2}{f_Z(z)S_n(C,h,z)},\;C\in
\CC,\;z\in H,\;h\in [h_n,\wth_n] \label{Ungz},\eeq and in
establishing the almost sure limit behaviour of the quantities
$U_n(C,h,z)$. Recall that $\sig^2(C,z):=\mathrm{Var}\po 1_C(Y)\mid
Z=z\pf$ and write \beq r(C,z):=\EEE\po 1_C(Y)\mid
Z=z\pf.\label{rCz}\eeq By $(Hf)$ together with Scheffé's lemma,
both $\sig^2(C,\cdot)$ and $r(C,\cdot)$ are equicontinuous
uniformly in $C\in \CC$, namely
\begin{align} \lim_{\dd\rar
0}\sup_{C\in\CC}\;\mathop{\sup_{z_1,z_2\in H}}_{\mid
z_1-z_2\mid\le \dd}\mid r(C,z_1)-r(C,z_2)\mid=0&,\label{unifr}\\
\lim_{\dd\rar 0}\sup_{C\in\CC}\;\mathop{\sup_{z_1,z_2\in H}}_{\mid
z_1-z_2\mid\le \dd}\mid
\sig^2(C,z_1)-\sig^2(C,z_2)\mid=0&\label{unifsig}.
\end{align}
\subsection{Step 1: an application of Theorem 2}
Recall that $\sig^2(C,z):=\mathrm{Var}\po 1_C(Y_i)\mid Z=z\pf$ and
that $r(C,z):=\PPP\po Y\in C\mid Z=z\pf$. In this first step we
prove that, given $\e>0$, we have
$(2\log(h^{-d}))^{-1}U_n(C,h,z)\le (1+\e)$ uniformly in $C,h,z$,
ultimately as $\nif$.
\begin{lem}\label{Sncv}
Under the assumptions of Theorem \ref{T2}, we have almost surely :
\begin{align} \limn&\mathop{\sup_{z\in H,\;C\in \CC,}}_{h_n\le h\le \wth_n}
\Mid \frac{S_n(C,h,z)}{nh^d\sig^2(C,z)\mmi K\mmi_{\lab,2}^2}-1\Mid=0,\label{poi1}\\
\limn&\mathop{\sup_{z\in H,\;C\in \CC,}}_{h_n\le h\le \wth_n} \Mid
\frac{X_n(C,h,z)}{\sqrt{2f_Z(z)\sig^2(C,z)\mmi
K\mmi_{\lab,2}^2nh^d\log(h^{-d})}}\Mid=1\label{poi2}. \end{align}
As a consequence we have
\begin{align}
\limn &\mathop{\sup_{z\in H,\;C\in \CC,}}_{h_n\le h\le \wth_n}
\Mid\frac{U_n(C,h,z)}{2\log
(h^{-d})}\Mid=1\;\;a.s.\label{pou}\end{align}
\end{lem}
\textbf{Proof}:\lb Note that (\ref{pou}) is a consequence of
(\ref{poi1}) and (\ref{poi2}). Set $L(\cdot)=K^2(\cdot)\mmi
K\mmi_{\lab,2}^{-2}$. To apply Theorem \ref{T1} we write $\po
1_C(Y)-r(C,z)\pf^2=1_C(Y)(1-2r(C,z))+r^2(C,z).$ Notice that, under
$(HG')$ and $(HV')$, the class $\CC$ and the sequence
$(h_n)\suite$ satisfy the conditions of Theorem \ref{T1} with
$p=\infty$. By Scheffé's lemma together with assumption $(Hf)$ and
(HG'), the two following collections of functions are uniformly
equicontinuous on $H$ : \beq \DD_1:=\ao
f_Z^{-1/2}(\cdot)(1-2r(C,\cdot)),\;C\in \CC \af,\;\;\DD_2:=\ao
f_Z^{-1/2}(\cdot)r^2(C,\cdot),\;C\in\CC\af.\eeq We can hence apply
Theorem \ref{T1} to the class $\CC$, with $\DD_1,\DD_2$ defined as
above, and with the kernel $L$ to obtain, with probability one,
\begin{align}
 &\limn \sup_{z\in H,\;C\in \CC,\;h_n\le h\le \wth_n} \frac{\mid
\wtW_n(C,h,z)\mid}{\sqrt{2nh^d\log(h^{-d})}}<\infty\label{const},
\end{align}
with \beq \nono\wtW_n(C,h,z):=f_Z(z)^{-1}\sliin \aooo\po
1_C(Y_i)-r(C,z)\pf^2L\poo \frac{Z_i-z}{h}\pff-\EEE\cooo\po
1_C(Y_i)-r(C,z)\pf^2L\poo \frac{Z_i-z}{h}\pff\cfff\afff.\eeq Now
write \beq \EEE\pooo \po 1_C(Y_i)-r(C,z)\pf^2 L\poo
\frac{Z_i-z}{h}\pff\pfff=: \wt{r}(C,h,z).\eeq By assumptions
$(HG')$, $(HV')$ and $(Hf)$ together with Scheffé's lemma, we can
infer that
\begin{align}
&\limn\mathop{\sup_{z\in H,\;C\in \CC,}}_{ h\in[h_n, \wth_n]}\Mid
\frac{\wt{r}(C,h,z)}{h^d f_Z(z)\sig^2(C,z)}-1\Mid=0\label{coco},\\
&\limn \mathop{\sup_{z\in H,\;C\in \CC,}}_{h\in [h_n,\wth_n]}\mid
m(C,h,z)-r(C,z)\mid=0,\label{moments}\\
&\limn \sup_{h\in
[h_n,\wth_n]}\frac{\log(h^{-d})}{nh^d}=0.\label{rapo}
\end{align}
Writing
\begin{align}
\nono &\poo S_n(C,h,z)-nf_Z(z)^{-1}\mmi
K\mmi_{\lab,2}^2\wt{r}(C,h,z)\pff-\wt{W}_n(C,h,z)\\\nono=&f_Z(z)^{-1}\coo
\po m^2(C,h,z)-r^2(C,z)\pf-2\po m(C,h,z)-r(C,z)\pf\cff\mmi
K\mmi^2_{\lab,2}\sliin L\poo \frac{Z_i-z}{h}\pff,
\end{align}
we conclude by Theorem \ref{T1} and (\ref{moments}) that
$$\limn\mathop{\sup_{z\in H,\;C\in\CC,}}_{h_n\le h\le \wth_n} \frac{\Mid
\poo S_n(C,h,z)-nf_Z(z)^{-1}\mmi
K\mmi_{\lab,2}^2\wt{r}(C,h,z)\pff-\wt{W}_n(C,h,z)\Mid}{\sqrt{2nh^d\log(h^{-d})}}=0$$
with probability one, from where we obtain with (\ref{const}) and
(\ref{coco}) that \begin{align}  \limn\mathop{\sup_{z\in H,\;C\in
\CC,}}_{h_n\le h\le \wth_n} \frac{\Mid S_n(C,h,z)-nh^d\mmi
K\mmi_{\lab,2}^2\sig^2(C,z)\Mid}{\sqrt{2nh^d\log(h^{-d})}}<\infty.\label{const1}
\end{align}
The proof of (\ref{poi1}) is now concluded, by (\ref{rapo}),
(\ref{const1}) and $(HG')$. Assertion (\ref{poi2}) can be proved
in a very similar way, taking care that the class $\DD:=\ao
f_Z(\cdot)^{-1/2}\sig(C,\cdot)^{-1}\af$ is uniformly
equicontinuous and bounded away from zero and infinity on $H$. We
omit details. $\Box$\lb
\subsection{Step 2: convex hull condition}
 The second step of our proof
of Theorem \ref{T2} is usually called the "convex hull condition".
\begin{lem}\label{ch}
With probability one, we have, for all large $n$ and for all $C\in
\CC,\;z\in H,\;h_n\le h\le \wth_n$, \beq \sharp \aoo i:\;
K\poo\frac{Z_i-z}{h}\pff \po 1_C(Y_i)-m(C,h,z)\pf>0\aff \in
\{1,2,\ldots,n-1\}\label{convexhull}.\eeq
\end{lem}
\textbf{Proof}: It is sufficient to prove that \beq \lin
\mathop{\inf_{z\in H,\;C\in\CC,}}_{h\in[h_n,\wth_n]}\PPP\pooo\pm
\po 1_C(Y)-m(C,h,z)\pf K\poo \frac{Z-z}{h}\pff>0\pfff>0
,\label{stp}\eeq and that the following class is
Glivenko-Cantelli: \beq\nono \AAA:=\aooo \aoo (y,\wtz)\in
\RRR^{d'}\times\RRR^d,\;\po 1_C(y)-m(C,h,z)\pf
K\poo\frac{\wtz-z}{h}\pff>0\aff,\;C\in \CC,\;h>0,\;z\in
H\afff.\eeq First note that $\AAA\subset \BB$, where \beq
\nono\BB:=\aooo \aoo (y,\wtz)\in \RRR^{d'}\times\RRR^d,\;\po
1_C(y)-a\pf K\poo\frac{\wtz-z}{h}\pff>0\aff,\;C\in
\CC,\;h>0,\;z\in \RRR^d,\;a\in \RRR\afff.\eeq By (HK1) and by
Lemma 2.6.18 in Van der Vaart and Wellner (1996), the two following classes of sets
are VC:
$$\BB_{\pm}:=\aoo \ao\wtz\in \RRR^d,\;\pm K\po h^{-d}(\wtz-z)\pf>0\af,\;z\in
\RRR^d,\;h>0\aff.$$ Moreover, as $\CC$ is a VC class of sets, we
straightforwardly deduce that the following class is also VC:
$$\MM_{\GG}:=\aoo \{z\in \chi,\;1_C(z)>a\},\;C\in \CC,\;a\in
\RRR\aff.$$ By a combination of points $(i)$ and $(ii)$ of Lemma
2.6.17 in  Van der Vaart and Wellner (1996), we conclude that $\BB$ is VC, which
entails that $\AAA$ is Glivenko-Cantelli. We now have to prove
(\ref{stp}). Define the following family of random variables
$$\HH_{\mh}:=\aoo\po 1_C(Y)-m(C,h,z)\pf
K\poo\frac{Z-z}{h}\pff,\;z\in H,\;C\in \CC,\;0< h\le \mh\aff.$$ By
the Cauchy-Schwarz inequality we have $\PPP(X>0)\geq
\EEE(X^2)^{-1}\EEE(X1_{X>0})^2$. Hence it is sufficient to prove
that, for $\mh$ small enough we have
\begin{align} &\inf_{X\in \HH_\mh} \EEE\poo
X1_{X>0}\pff=\frac{1}{2}\inf_{X\in \HH_\mh} \EEE\poo \mid X \mid
\pff>0,\label{iinf}\\
&\sup_{X\in \HH_\mh} \EEE\poo X^2 \pff<\infty.\label{ssup}
\end{align}
Note that the equality appearing in (\ref{iinf}) is a consequence
of $\EEE(X)=0$ for each $X\in \HH_\mh$. By $(HG')$, $(Hf)$ and
(\ref{unifr}), routine analysis shows that, for $\mh$ small
enough, both (\ref{ssup}) and the following assertion are true:
\beq \inf_{X\in \HH_h}\EEE\po X^2\pf>\frac{1}{2}\inf_{z\in
H,\;C\in \CC} \sig^2(C,z)f_Z(z)\mmi
K\mmi_{\lab,2}=:\alp_0>0.\label{msn}\eeq Now, as $\HH_\mh$ is
uniformly bounded by some constant $M>0$ we get that $\alp_0\le M
\EEE\po\mid X\mid\pf$ for all $X\in \HH_\mh$, and hence
(\ref{iinf}) is proved. This concludes the proof of Lemma
\ref{ch}. $\Box$
\subsection{Step 3: end of the proof of Theorem \ref{T2}}
Lemma \ref{ch} ensures us (see, e.g., Owen (2001), p. 219) that
almost surely, for all large $n$ and for each $z\in
H,\;C\in\CC,\;h_n\le h\le \wth_n $, the maximum value in $\RR_n\po
m(C,h,z),C,h,z\pf$ is obtained by choosing the following weights
(recall (\ref{RRnthtgz})):
\beq
p_i(C,h,z):=
\frac{1}{n}\frac{1}{1+\lab_n(C,h,z)w_{i,n}(C,h,z)},\label{pigz}\eeq
where $\lab_n(C,h,z)$ is the unique solution of \beq \sliin
\frac{w_{i,n}(C,h,z)}{1+\lab_n(C,h,z)w_{i,n}(C,h,z)}=0.\label{equationlab}\eeq
Our next lemma gives an asymptotic control of $$
\mathop{\sup_{C\in \CC,z\in H,}}_{h_n\le h\le
\wth_n}\mid\lab_n(C,h,z)\mid.$$ It is largely inspired by Lemma 1
in Chen \textit{et al.} (2003).
\begin{lem}\label{lemlab}
Under the assumptions of Theorem \ref{T1} we have almost surely:
\beq \mathop{\sup_{C\in \CC,\;z\in H,}}_{h_n\le h\le \wth_n}
\sqrt{\frac{n h^d}{\log(h^{-d})}}\mid\lab_n(C,h,z)\mid = O(1).\eeq
\end{lem}
\textbf{Proof :} Following the proof of Owen (2001), p.
220, Lemma \ref{lemlab} will be proved if we check the following
three conditions:
\begin{align}
\max_{1\le i\le n}\;\mathop{\sup_{z\in H,\;C\in\CC,}}_{h_n\le h\le
\wth_n}\sqrt{\frac{\log(h^{-d})}{nh^d}}\mid w_{i,n}(C,h,z)\mid=&
o_{a.s.}(1)\label{c1},\\
\mathop{\sup_{z\in H,\;C\in\CC,}}_{h_n\le h\le
\wth_n} \frac{\mid X_n(C,h,z)\mid}{\sqrt{nh^d\log(h^{-d})}}= &O_{a.s.}(1),\label{c2}\\
\lin \mathop{\inf_{z\in H,\;C\in\CC,}}_{h_n\le h\le \wth_n}
\frac{S_n(C,h,z)}{nh^d}>&0\;a.s.\label{c3}
\end{align}
As each $w_{i,n}(C,h,z)$ is almost surely bounded by $2\mmi
K\mmi_{\RRR^d}$, and by (\ref{rapo}), condition (\ref{c1}) is
readily satisfied. Now note that condition (\ref{c2}) is a
straightforward consequence of Theorem \ref{T1}, and that
(\ref{c3}) is a consequence of both Lemma \ref{Sncv} and $(HG')$.
The remainder of the proof of Lemma \ref{lemlab} is done by
following Owen (2001), p. 220.$\Box$\vk Now set
$$V_{i,n}(C,h,z):=\lab_n(C,h,z)w_{i,n}(C,h,z).$$
By Lemma \ref{lemlab} and assertion (\ref{c1}) we have \beq \limn
\max_{1\le i\le n} \mathop{\sup_{z\in H,\;C\in\CC,}}_{h_n\le h\le
\wth_n} \mid V_{i,n}(C,h,z)\mid
= 0\text{ a.s. },\label{trup}\eeq which entails, almost surely,
for all large $n$ and for each $z\in H,\;C\in \CC,\;h\in
[h_n,\mh_n]$:
\begin{align}
\nono 0= &\sliin \frac{w_{i,n}(C,h,z)} {1+V_{i,n}(C,h,z)}\\
\nono =&\sliin w_{i,n}(C,h,z)\poo
1-V_{i,n}(C,h,z)+V^2_{i,n}(C,h,z)/(1+V_{i,n}(C,h,z))\pff\\
\nono=& X_n(C,h,z)-f_Z(z)S_n(C,h,z)\lab_n(C,h,z)+\sliin
\frac{w_{i,n}(C,h,z)V_{i,n}^2(C,h,z)}{1+V_{i,n}(C,h,z)}\\
=&X_n(C,h,z)-f_Z(z)S_n(C,h,z)\lab_n(C,h,z)+\sliin
\frac{w_{i,n}^3(C,h,z)}{1+V_{i,n}(C,h,z)}\lab^2_n(C,h,z).\label{rjv}
\end{align}
From (\ref{c1}), (\ref{c2}) and (\ref{trup}), we conclude that
there exists a random sequence $\e_n$ such that, almost surely, we
have $\e_n \rar 0$ and
\begin{align}
 \nono\sliin
\frac{w_{i,n}^3(C,h,z)}{1+V_{i,n}(C,h,z)}\lab^2_n(C,h,z)\le
&X_n(C,h,z) \max_{i}w^2_{i,n}(C,h,z) \\
\nono &\;\times \poo\min_{1\le i\le n} \mid
1+V_{i,n}(C,h,z)\mid\pff^{-1}\lab^2_n(C,h,z)\\
\le &\e_n \sqrt{nh^d\log(h^{-d})},
\end{align}
uniformly in $C\in \CC,z\in H,\;h\in [h_n,\wth_n]$. Hence,
dividing the right hand side of $(\ref{rjv})$ by $S_n(C,h,z)$,
recalling (\ref{c3}) and (\ref{rapo}), we obtain with probability
one that
\beq
\lab_n(C,h,z)=\frac{X_n(C,h,z)}{f_Z(z)S_n(C,h,z)}+\beta_n(C,h,z),
\eeq with $\beta_n(C,h,z)\le M\e_n\sqrt{\log(h^{-d})/nh^d}$
uniformly in $C\in \CC,\;z\in H$ and $h\in [h_n,\wth_n]$, for some
almost surely finite random variable $M$. We can now conclude that
(recall (\ref{Ungz})) \beq\limn \mathop{\sup_{z\in H,\;C\in
\CC,}}_{h\in [h_n,\wth_n]}\Mid \frac{-2\log\poo\RR_n\po
g,z,m(C,h,z)\pf\pff}{U_n(C,h,z)}-1\Mid=0\label{dajou1},\eeq by
reasoning as in Owen (2001), p. 221. The proof of Theorem
\ref{T2} is then concluded by (\ref{pou}). $\Box$

\end{document}